\RequirePackage{ifpdf}
\ifpdf 
\documentclass[pdftex]{sigma}
\else
\documentclass{sigma}
\fi

\usepackage[mathscr]{eucal}

\begin{document}

\allowdisplaybreaks

\renewcommand{\PaperNumber}{028}

\FirstPageHeading

\renewcommand{\thefootnote}{$\star$}

\ShortArticleName{Bethe Ansatz for the Ruijsenaars Model of
$BC_1$-Type}

\ArticleName{Bethe Ansatz for the Ruijsenaars Model of
$\boldsymbol{BC_1}$-Type\footnote{This paper is a contribution to
the Vadim Kuznetsov Memorial Issue `Integrable Systems and Related
Topics'. The full collection is available at
\href{http://www.emis.de/journals/SIGMA/kuznetsov.html}{http://www.emis.de/journals/SIGMA/kuznetsov.html}}}

\Author{Oleg CHALYKH}

\AuthorNameForHeading{O. Chalykh}

\Address{School of Mathematics, University of Leeds, Leeds LS2
9JT, United Kingdom}
\Email{\href{mailto:o.chalykh@leeds.ac.uk}{o.chalykh@leeds.ac.uk}}

\ArticleDates{Received December 14, 2006, in f\/inal form February
06, 2007; Published online February 22, 2007}

\Abstract{We consider one-dimensional elliptic Ruijsenaars model
of type $BC_1$. It is given by a three-term dif\/ference
Schr\"odinger operator $L$ containing $8$ coupling constants. We
show that when all coupling  constants are integers, $L$ has
meromorphic eigenfunctions expressed by a variant of Bethe ansatz.
This result generalizes the Bethe ansatz formulas known in the
$A_1$-case.}

\Keywords{Heun equation; three-term dif\/ference operator; Bloch
eigenfunction; spectral curve}

\Classification{33E30; 81U15}

\def\c{\mathbb{C}}
\def\e{\boldsymbol{e}}
\def\a{\mathbb{A}}
\def\k{\mathbb{K}}
\def\Z{\mathbb{Z}}
\def\Q{\mathbb{Q}}
\def\N{\mathbb{N}}
\def\P{\mathbb{P}}
\def\D{\mathcal{D}}
\def\O{\mathcal{O}}
\def\h{\mathfrak{h}}
\def\g{\mathfrak{g}}

\begin{flushright}
\textit{Dedicated to the memory of Vadim Kuznetsov}
\end{flushright}

\section{Introduction}

The quantum Ruijsenaars model \cite{R} in the simplest two-body
case reduces to the following dif\/ference operator acting on
functions of one variable:
\begin{equation}\label{ru}
L=\frac{\sigma(z-2\gamma m)}{\sigma(z)}T^{2\gamma}+
\frac{\sigma(z+2\gamma m)}{\sigma(z)}T^{-2\gamma},
\end{equation}
where $\sigma(z)$ is the Weierstrass $\sigma$-function, $m$ is the
coupling parameter, and $T^\gamma$ stands for the shift operator
acting by $(T^{\gamma}f)(z)=f(z+\gamma)$. This operator, which
f\/irst appeared in E.~Sklyanin's work~\cite{S,S1}, can be viewed
as a dif\/ference version of the Lam\'e operator
$-d^2/dz^2+m(m+1)\wp(z)$. It was observed by Krichever--Zabrodin
\cite{KZ} and by Felder--Varchenko \cite{FV1, FV2}, that in the
special case of integer coupling parameter the operator \eqref{ru}
shares many features with the Lam\'e operator. In particular, when
$m\in\Z_+$ they both have meromorphic Bloch eigenfunctions which
can be given explicitly by a suitable Bethe ansatz. This
ref\/lects the well-known fact that the Lam\'e operator is {\it
finite-gap} for integer $m$ (see e.g.~\cite{DKN} for a survey of
the f\/inite-gap theory).

The Lam\'e operator has the following generalization closely
related to the Heun's equation:
\begin{equation}\label{h}
H=-d^2/dz^2+\sum_{p=0}^3 g_p(g_p+1)\wp(z+\omega_p),
\end{equation}
where $\omega_p$ are the half-periods of $\wp(z)$. It can be
viewed as a $BC_1$-generalization of the Lam\'e operator. Again,
for integer coupling parameters $g_p$ this operator is
f\/inite-gap, as was discovered by Treibich--Verdier \cite{TV},
see also \cite{Sm,T} for the detailed study of \eqref{h}.

The operator \eqref{h} has a multivariable generalization known as
the Inozemtsev model \cite{I}. A~relativistic version of the
Inozemtsev model ($\equiv$ $BC_n$ version of the Ruijsenaars
model) was suggested by J.F.~van Diejen \cite{vD1, vD2}, see also
\cite{KH1,KH2}. In the simplest one-variable case it takes the
form of a three-term dif\/ference operator
\begin{gather}\label{l}
L=a(z)T^{2\gamma}+b(z)T^{-2\gamma}+c(z), \intertext{where}
\label{a}
a(z)=\prod_{p=0}^3\frac{\sigma_{p}(z-\mu_p)\sigma_{p}(z+\gamma-\mu'_p)}
{\sigma_{p}(z)\sigma_{p}(z+\gamma)},\qquad b(z)=a(-z),
\end{gather}
$c(z)$ is given explicitly in \eqref{c} below, and the notations
are explained at the beginning of the next section.

Therefore, the operator \eqref{l} should be viewed both as a
dif\/ference analogue of \eqref{h} and a $BC_1$-version of the
Ruijsenaars model \eqref{ru}. In the trigonometric limit it
coincides with the Askew--Wilson dif\/ference operator \cite{AW}.
(Notice that \eqref{l} contains eight parameters $\mu_p$,
$\mu_p'$, compared to the four in the Askew--Wilson operator and
in \eqref{h}.) Therefore, it is natural to expect that \eqref{l}
and \eqref{ru} should have similar properties. This is indeed the
case, as we will demonstrate below. Our main result says that in
the case of integer coupling parameters
\begin{equation}\label{cc}
\mu_p=2\gamma m_p,\qquad \mu'_p=2\gamma m'_p,\qquad m_p,
m_p'\in\Z_+\qquad (p=0,\dots ,3),
\end{equation}
the operator \eqref{l} has meromorphic Bloch eigenfunctions which
can be given explicitly via a version of Bethe ansatz. Note that
our derivation of the Bethe ansatz equations is very elementary:
it only uses some simple facts about the operator \eqref{l}.

\section[Ruijsenaars operator of type $BC_1$]{Ruijsenaars operator of type $\boldsymbol{BC_1}$}

\subsection{Preliminaries}
Let $\sigma(z)=\sigma(z;2\omega_1,2\omega_2)$ denote the
Weierstrass $\sigma$-function with the half-periods $\omega_1$,
$\omega_2$. Recall that $\sigma(z)$ is an entire odd function on
the complex plane quasiperiodic with respect to $2\omega_1$,
$2\omega_2$. It will be convenient to introduce the third half
period as $\omega_3=-\omega_1-\omega_2$. One has
\begin{equation*}\label{qs}
\sigma(z+2\omega_s)=-\sigma(z)e^{2\eta_s(z+\omega_s)},\qquad
s=1,2,3,
\end{equation*}
with $\eta_s=\zeta(\omega_s)$, where
$\zeta(z)=\sigma'(z)/\sigma(z)$ denotes the Weierstrass
$\zeta$-function. Clearly, $\eta_1+\eta_2+\eta_3=0$. There is a
relation between $\eta_s$ and half-periods as follows:
\begin{equation*}
\eta_1\omega_2-\eta_2\omega_1=\pi i/2\,.
\end{equation*}
It is known that $\sigma(z)$ has simple zeros at points of the
period lattice $2\Gamma$, where $\Gamma=\omega_1\Z+\omega_2\Z$.
Likewise, $\zeta(z)$ has simple poles with residue $1$ at those
points. Let us introduce the shifted versions of $\sigma$ as
follows:
\begin{equation*}\label{ss}
\sigma_r(z)=e^{-\eta_rz}\sigma(z+\omega_r)/\sigma(\omega_r),\qquad
r=1,2,3.
\end{equation*}
These are even functions: $\sigma_r(-z)=\sigma_r(z)$. The
quasiperiodicity properties of $\sigma_r(z)$ look as follows:
\begin{equation*}\label{qss}
\sigma_r(z+2\omega_s)=(-1)^{\delta_{r,s}}\sigma_r(z)e^{2\eta_s(z+\omega_s)},\qquad
r, s=1,2,3.
\end{equation*}
Below we will use the convention that $\sigma_0(z)=\sigma(z)$ and
$\omega_0=\eta_0=0$. Note that
\begin{equation*}\label{s2}
\sigma(2z)=2\sigma_0(z)\sigma_1(z)\sigma_2(z)\sigma_3(z).
\end{equation*}

Let us remark on quasiperiodicity of the coef\/f\/icients $a(z)$
and $b(z)=a(-z)$ of the operator~\eqref{l}. Let us use the
subscript to indicate the dependence of $a(z)=a_\mu(z)$ and
$b(z)=b_\mu(z)$ on the parameters $\mu=\{\mu_p, \mu_p'\}$. Using
the translation properties of $\sigma_p$ one checks directly that
$a(z)$, $b(z)$ given by \eqref{a} have the following covariance
with respect to the shift by a half-period $\omega_r$ ($r=0,\dots,
3$):
\begin{equation}\label{co}
a_\mu(z+\omega_r)=a_{\widetilde\mu}(z)e^{\eta_r\sum\limits_{p=0}^3(\mu_p+\mu'_p)},\qquad
b_\mu(z+\omega_r)=b_{\widetilde\mu}(z)e^{-\eta_r\sum\limits_{p=0}^3(\mu_p+\mu'_p)},
\end{equation}
with $\widetilde\mu_p=\mu_{\pi_r(p)}$ and
$\widetilde\mu'_p=\mu'_{\pi_r(p)}$, where $\pi_r$ is one of the
following permutations:
\begin{equation}\label{pi}
\pi_0=\mathtt{id},\qquad \pi_1=(01)(23),\qquad
\pi_2=(02)(13),\qquad \pi_3=(03)(12).
\end{equation}
Note that these permutations form an Abelian subgroup of $S_4$,
and $\pi_p(0)=p$ for $p=0,\dots, 3$. It follows that
$\pi_r\circ\pi_p=\pi_q$ whenever $q=\pi_r(p)$.

\subsection[Ruijsenaars operator of type ${BC}_1$]{Ruijsenaars operator of type $\boldsymbol{{BC}_1}$}
Let $L$ be the operator \eqref{l}--\eqref{a} with the
coef\/f\/icient $c(z)$ given by
\begin{equation}\label{c}
c(z)=\sum_{p=0}^3 c_p(\zeta_p(z+\gamma)-\zeta_p(z-\gamma)),\qquad
\zeta_p(z)=\frac{\sigma_p'(z)}{\sigma_p(z)}=-\eta_p+\zeta(z+\omega_p),
\end{equation}
where $c_p$ looks as follows:
\begin{equation}\label{cp}
c_p=-\frac{2}{\sigma(2\gamma)}\prod_{s=0}^3\sigma_s(\gamma+\mu_{\pi_p(s)})\sigma_s(\mu'_{\pi_p(s)}).
\end{equation}
Here the permutations $\pi_p$ are the same as in \eqref{pi}.

\begin{remark}
$L$ is a one-dimensional counterpart of its two-variable version
introduced in \cite{vD1}, see also \cite{vD2,KH1,KH2}. Their
precise relation is as follows: the $BC_2$ version in \cite{vD1}
involves one more parameter $\mu$ attached to the roots $e_1\pm
e_2$, and it decouples when $\mu=0$. Namely, let $D$ be as
in~\cite{vD1}, formulas (3.11), (3.19)--(3.20), (3.23)--(3.24).
Then we have, for an appropriate constant $\alpha$, that
$\lim\limits_{\mu\to 0} (D- \alpha \mu^{-1})=L_1+L_2$, where
$L_1$, $L_2$ act in $z=z_1$ and $z_2$, respectively, and are given
by \eqref{l}--\eqref{a} and \eqref{c}--\eqref{cp} above. Observe
that in the special case $\mu_0=2\gamma m$,
$\mu_0'=\mu_p=\mu_p'=0$ ($p=1,2,3$) the operator $L$ reduces to
\eqref{ru}.
\end{remark}

\subsection[Symmetries of $L$]{Symmetries of $\boldsymbol{L}$}
It is obvious from the formula \eqref{c} for $c(z)$ that
$c(-z)=c(z)$. As a result, $L$ is invariant under
$(z\leftrightarrow -z)$. Next, a direct computation shows that the
function $c(z)=c_\mu(z)$ is covariant under the shifts by
half-periods, namely:
\begin{equation}\label{cco}
c_\mu(z+\omega_r)=c_{\widetilde \mu}(z),\qquad\text{where}\quad
\widetilde\mu_p=\mu_{\pi_r(p)},\quad
\widetilde\mu'_p=\mu'_{\pi_r(p)}.
\end{equation}
Let us write $L=L_\mu$ to indicate dependence on $\mu$. Combining
\eqref{co} and \eqref{cco}, we obtain that
\begin{equation*}\label{col}
T^{\omega_r}\circ L_\mu\circ T^{-\omega_r}=e^{-\lambda_r z}\circ
L_{\widetilde\mu}\circ e^{\lambda_r z},\qquad
\lambda_r=\eta_r({2\gamma})^{-1}\sum_{p=0}^3 (\mu_p+\mu'_p),
\end{equation*}
where $\widetilde\mu=\pi_r(\mu)$ is the same as in \eqref{co},
\eqref{cco}. This implies that for any $\omega\in\Gamma$
\begin{equation}\label{col2}
T^{\omega}\circ L_\mu\circ T^{-\omega}=e^{-\lambda(\omega) z}\circ
L_{\widetilde\mu}\circ e^{\lambda(\omega)z},
\end{equation}
where $\lambda(\omega):=n_1\lambda_1+n_2\lambda_2$ if
$\omega=n_1\omega_1+n_2\omega_2$, and $\widetilde \mu$ in the
right-hand side is def\/ined as $\widetilde \mu=\pi_s(\mu)$ if
$\omega\equiv \omega_s \mod 2\Gamma$.

\section{Bethe ansatz}

Providing the coupling constants satisfy \eqref{cc}, put
$m=\sum\limits_{p=0}^3(m_p+m'_p)$ and consider the following
function $\psi(z)$ depending on the parameters $t_1,\dots, t_m,
k\in\c$:
\begin{equation}\label{psi}
\psi(z)=e^{kz}\prod_{i=1}^m\sigma(z+t_i).
\end{equation}
Let us impose $m$ relations onto these parameters as follows:
\begin{gather}\label{ba}
\psi(\omega_s+2j\gamma)= \psi(\omega_s-2j\gamma)e^{4j\gamma
m\eta_s} \qquad (j=1,\dots, m_s),\\ \label{ba'}
\psi(\omega_s+(2j-1)\gamma)=
\psi(\omega_s-(2j-1)\gamma)e^{(4j-2)\gamma m\eta_s} \qquad
(j=1,\dots, m'_s).
\end{gather}
(Here $s=0,\dots, 3$.) We will refer to \eqref{ba}--\eqref{ba'} as
the {\it Bethe ansatz equations}, or simply the Bethe equations.
Explicitly, they look as follows:
\begin{gather}\label{bae}
e^{4j\gamma
m\eta_s}\prod_{i=1}^m\frac{\sigma(t_i+\omega_s-2j\gamma)}{\sigma(t_i+\omega_s+2j\gamma)}=
e^{4j\gamma k}\qquad (j=1,\dots, m_s),\\ \label{bae'}
e^{(4j-2)\gamma
m\eta_s}\prod_{i=1}^m\frac{\sigma(t_i+\omega_s-(2j-1)\gamma)}{\sigma(t_i+\omega_s+(2j-1)\gamma)}=
e^{(4j-2)\gamma k}  \qquad (j=1,\dots, m'_s).
\end{gather}

Now we can formulate the main result of this paper.

\begin{theorem}\label{mt} Suppose the parameters $t_1,\dots, t_m, k$ satisfy the
Bethe equations \eqref{ba}--\eqref{ba'} and the conditions
$t_i+t_j\notin 2\omega_1\Z+2\omega_2\Z$ for $1\le i\ne j\le m$.
Then the corresponding function $\psi(z)$ \eqref{psi} is an
eigenfunction of the operator \eqref{l}--\eqref{cc}.
\end{theorem}
The proof will be given in the next section.

\begin{remark}
To compute the corresponding eigenvalue, one evaluates the
expression $L\psi/\psi$ at any suitable point $z$. For instance, a
convenient choice is $z=2\gamma m_0$ (provided $m_0>0$), because
then the f\/irst term in $L\psi$ vanishes.
\end{remark}

\begin{remark} If some of the coupling parameters $m_p$, $m'_p$ vanish,
then the corresponding sets of the Bethe equations are not present
in \eqref{ba}--\eqref{ba'}. For example, in the case when the only
nonzero parameter is $m_0=m$, the Bethe equations take the form:
\begin{equation*}\label{ba0}
\psi(2j\gamma)=\psi(-2j\gamma) \qquad (j=1,\dots, m).
\end{equation*}
In that form (seemingly dif\/ferent from \cite{KZ,FV1}) they
appeared in \cite{Z}.
\end{remark}

\subsection{Invariant subspaces}

The idea of the proof of the theorem is that applying $L$ to
$\psi$ will not destroy the conditions \eqref{ba}--\eqref{ba'},
cf.~\cite{C,CEO}. We begin with two elementary results about a
three-term dif\/ference operator with meromorphic
coef\/f\/icients:
\begin{equation*}\label{d}
D=a(z)T^{2\gamma}+b(z)T^{-2\gamma}+c(z).
\end{equation*}
Suppose that $a$, $b$, $c$ are regular at $z\in 2\gamma\Z$, apart
from $z=0$ where $a$, $b$ have simple poles. Furthermore, suppose
that
\begin{gather}\label{re}
\mathrm{res}_{z=0}(a+b)=0 \intertext{and that for some $m\in\Z_+$
the following is true:} \label{c1} a(2\gamma m)=0,\qquad
a(2j\gamma)=b(-2j\gamma),\qquad c(2j\gamma)=c(-2j\gamma) \qquad
(j=\pm 1, \dots, \pm m).
\end{gather}

\begin{lemma}[cf.~\cite{C}, Lemma 2.2]\label{le1} Let $D$ be as above and define
$Q_m$ as the space of meromorphic functions $f(z)$ which are
regular at all points $z\in 2\gamma\Z$ and satisfy the conditions
$f(2j\gamma)=f(-2j\gamma)$ for all $j=1,\dots, m$. Then
$D(Q_m)\subseteq Q_m$.
\end{lemma}

\begin{proof} For $D'=aT^{2\gamma}+bT^{-2\gamma}$ this is precisely Lemma 2.2 from \cite{C},
thus $D'(Q_m)\subseteq Q_m$. On the other hand, the conditions on
$c$ in \eqref{c1} imply trivially that $cQ_m\subseteq Q_m$.
\end{proof}

\begin{corollary}\label{cor1}
Suppose that instead of \eqref{re}, \eqref{c1} we know that $D$ is
invariant under $(z\leftrightarrow -z)$ and that $a(2\gamma m)=0$.
Then $D(Q_m)\subseteq Q_m$.
\end{corollary}

\begin{proof}
Indeed, in that case we know that $b(z)=a(-z)$ and $c(z)=c(-z)$
identically. This implies the conditions \eqref{re}--\eqref{c1}.
\end{proof}

For the next lemma, we assume that: (1) $a$ is regular at $z\in
\gamma+2\gamma\Z$\, apart from a simple pole at $z=-\gamma$; (2)
$b$ is regular at $z\in \gamma+2\gamma\Z$\, apart from a simple
pole at $z=\gamma$; (3) $c$ is regular at $z\in \gamma+2\gamma\Z$
apart from simple poles at $z=\pm \gamma$. Also, suppose that
\begin{gather}\label{re2}
\mathrm{res}_{z=-\gamma}(a+c)=\mathrm{res}_{z=\gamma}(b+c)=0,\qquad
\mathrm{res}_{z=-\gamma}a=\mathrm{res}_{z=\gamma}b,\\\label{re3}
(a+b+c)|_{z=-\gamma}=(a+b+c)|_{z=\gamma}.
\end{gather}
(The last condition makes sense because \eqref{re2} implies that
$a+b+c$ is regular at $z=\pm \gamma$.) In addition to that ,
assume that for some $m\in\Z_+$ the following is true:
\begin{gather}
\label{c2} a((2m-1)\gamma)=0,\qquad
a((2j-1)\gamma)=b((-2j+1)\gamma)\quad \text{for} \ \ j=\pm
1,\dots, \pm
m,\\
\label{c3} c((2j+1)\gamma)=c(-(2j+1)\gamma)\quad \text{for }
j=1,\dots, m-1.
\end{gather}

\begin{lemma}[cf. \cite{C}, Lemma 2.3]\label{le2} For $D$ as above, def\/ine
$Q'_m$ as the space of meromorphic functions $f(z)$ which are
regular at $z\in \gamma+2\gamma\Z$\, and satisfy the conditions
$f((2j-1)\gamma)=f((-2j+1)\gamma)$ for $j=1,\dots, m$. Then
$D(Q'_m)\subseteq Q'_m$.
\end{lemma}

\begin{proof} For $D'=a(T^{2\gamma}-1)+b(T^{-2\gamma}-1)$ the proof of the inclusion
$D'(Q'_m)\subset Q'_m$ follows the proof of Lemma 2.3 in \cite{C}.
On the other hand, the conditions on $a$, $b$, $c$ imply that
$c':=c+a+b$ belongs to $Q'_m$. Thus, $c'Q'_m\subseteq Q'_m$.
Therefore, the operator $D=D'+(a+b+c)$ preserves~$Q'_m$.
\end{proof}

\begin{corollary}\label{cor2} The lemma above remains valid after replacing
\eqref{re3}--\eqref{c3} by the invariance of $D$ under
$(z\leftrightarrow -z)$ and the condition that
$a((2m-1)\gamma)=0$.
\end{corollary}

\begin{proof} Indeed, the conditions \eqref{re3}--\eqref{c3} follow easily
from the fact that $b(z)=a(-z)$ and $c(z)=c(-z)$.
\end{proof}

Let us apply these facts to the Ruijsenaars operator~\eqref{l}
with integer coupling parameters~\eqref{co}. Below we always
assume that the step $\gamma$ is irrational, i.e.\ $\gamma\notin
\mathbb Q\otimes_{\Z} \Gamma=\mathbb Q\omega_1+\mathbb Q\omega_2$.
We proceed by def\/ining $Q$ as the space of entire functions
$\psi(z)$ satisfying the following conditions for every $\omega\in
\omega_s+2\Gamma$ ($s=0,\dots, 3$):
\begin{gather}\label{bal}
\psi(\omega+2j\gamma)= \psi(\omega-2j\gamma)e^{4j\gamma
m\eta(\omega)} \qquad (j=1,\dots, m_s),\\ \label{bal'}
\psi(\omega+(2j-1)\gamma)=
\psi(\omega-(2j-1)\gamma)e^{(4j-2)\gamma m\eta(\omega)} \qquad
(j=1,\dots, m'_s).
\end{gather}
Here $m$ stands as before for $m=\sum\limits_{p=0}^3 (m_p+m'_p)$,
and the constant $\eta(\omega)$ is def\/ined for
$\omega=n_1\omega_1+n_2\omega_2$ as
$\eta(\omega)=n_1\eta_1+n_2\eta_2$.

\begin{proposition}\label{in} For integer coupling parameters \eqref{co} the Ruijsenaars
operator \eqref{l} preserves the space $Q$ of entire functions
with the properties \eqref{bal}--\eqref{bal'}:\ $L(Q)\subseteq Q$.
\end{proposition}

\begin{proof} First, by applying Corollaries \ref{cor1}, \ref{cor2} to the
Ruijsenaars operator, we obtain that $L$ preserves the spaces
$Q_{m_0}$ and $Q'_{m'_0}$ (in the notations of Lemmas \eqref{le1},
\eqref{le2}). Note that in doing so, we only have to check the
vanishing of $a(z)$ as required in Corollaries \ref{cor1},
\ref{cor2}, and the conditions on the residues \eqref{re2}. This
is where the formula \eqref{cp} becomes crucial. Finally, in order
to show that $L$ preserves similar conditions at other points
$\omega\in\Gamma$, one applies the formula \eqref{col2}.
\end{proof}

\medskip
Next, given $\alpha_1,\alpha_2\in\c$ and $m\in\Z$, let us write
$\mathcal F^{\alpha_1,\alpha_2}_m$ for the space of meromorphic
functions~$\psi(z)$ having the following quasiperiodicity
properties:
\begin{equation*}
\psi(z+2\omega_s)=e^{m\eta_sz+\alpha_s}\psi(z),\qquad s=1,2.
\end{equation*}
Entire functions in $\mathcal F^{\alpha_1,\alpha_2}_m$ ($m>0$) are
known as {\it theta-functions of order} $m$ (with
characteristics), each of them being a constant multiple of
\eqref{psi}, for appropriate $t_1,\dots, t_m, k$.

Now, a simple check shows that in the case \eqref{co} the
Ruijsenaars operator \eqref{l} preserves these spaces
corresponding to $m=\sum\limits_{p=0}^3(m_p+m'_p)$:
\begin{equation*}\label{qp}
L(\mathcal F^{\alpha_1,\alpha_2}_m)\subseteq \mathcal
F^{\alpha_1,\alpha_2}_m,\qquad\forall\,\alpha_1, \alpha_2.
\end{equation*}
Combining this with Proposition \ref{in}, we conclude that $L$
preserves the space of theta-functions of order $m$ satisfying the
conditions \eqref{bal}--\eqref{bal'}:
\begin{equation}\label{pr}
L(\mathcal F^{\alpha_1,\alpha_2}_m\cap Q)\subseteq \mathcal
F^{\alpha_1,\alpha_2}_m\cap Q\,,\quad\forall\,\alpha_1,
\alpha_2\,.
\end{equation}

\begin{proof}[Proof of the Theorem \ref{mt}] Take a solution
$(t_1,\dots, t_m, k)$ to the Bethe equations and the corresponding
function $\psi$ \eqref{psi}. Clearly, $\psi$ belongs to $\mathcal
F^{\alpha_1,\alpha_2}_m$ for some $\alpha_1, \alpha_2$. The Bethe
equations give the conditions \eqref{bal}--\eqref{bal'} only for
$\omega=\omega_s$, but the rest follows from the translation
properties of $\psi$. Thus, $\psi$ belongs to the space $\mathcal
F^{\alpha_1,\alpha_2}_m\cap Q$. By \eqref{pr},
$\widetilde\psi:=L\psi$ also belongs to this space. Now we use the
following fact (whose proof will be given below):
\begin{lemma}\label{l2} For any two functions $\psi, \widetilde\psi\in F^{\alpha_1,\alpha_2}_m\cap
Q$, their quotient $\widetilde\psi/\psi$ is an even elliptic
function, i.e.\ it belongs to $\c(\wp(z))$.
\end{lemma}
By this lemma, if $\widetilde\psi/\psi$ is not a constant, then
its poles must be invariant under $z\mapsto -z$, thus there exist
at least two of $t_1, \dots, t_m$ such that their sum belongs to
$2\Gamma$.
\end{proof}

\begin{proof}[Proof of the lemma]
Take any two functions $\psi, \widetilde\psi$ in $\mathcal
F^{\alpha_1,\alpha_2}_m\cap Q$ and put $f:=\widetilde\psi/\psi$.
Note that $f$ is an elliptic function of degree $\le m$ (because
its denominator and numerator have $m$ zeros in the fundamental
region). Let us label $m$ pairs of points $\omega_s\pm 2j\gamma$,
$\omega_s\pm (2j-1)\gamma$ as $P^\pm_l$ with $l=1,\dots, m$, then
the properties of $\psi$, $\widetilde\psi$ imply that $f$
satisf\/ies the conditions
\begin{equation}\label{balf}
f(P^+_l)=f(P^-_l),\qquad l=1,\dots, m.
\end{equation}
We may assume that $f$ is regular in at least one of the
half-periods $\omega_s$, otherwise switch to
$1/f=\psi/\widetilde\psi$. Let us anti-symmetrize $f$ to get
$g(z):=f(z)-f(-z)$, which will be odd elliptic, of degree $\le
2m$. It is clear that $g$ also satisf\/ies the conditions
\eqref{balf}. At the same time, it is anti-symmetric under any of
the transformations $z\mapsto 2\omega_s-z$ ($s=0,\dots, 3$).
Altogether this implies that $g$ must vanish at each of the $2m$
points $P^\pm_l$. Finally, it must vanish at one of the
half-periods (where $f$ was regular). So $g$ has $\ge 2m+1>\deg
(g)$ zeros, hence $g=0$, $f(z)\equiv f(-z)$, and we are done.

The above argument, however, would not work if one or both of the
functions $\psi$, $\widetilde \psi$ vanish at some of the points
$P^\pm_l$. Indeed, then we cannot claim that $f$ is regular at
those points, so some of the conditions \eqref{balf} would not
hold for $f$. In that case, we can argue as follows. Let
$\psi_\lambda$ denote the linear combination
$\psi_\lambda=\psi+\lambda\widetilde\psi$. Then $\psi_\lambda$,
$\psi_\mu$ for generic $\lambda$, $\mu$ will have zero of the same
multiplicity at any given point $P^\pm_l$. Thus, choosing
$\lambda$, $\mu$ appropriately, we can always achieve that
$\psi_\lambda/\psi_\mu\ne 0,\infty$ at every of these $2m$ points.
Let $r$ be the number of those pairs $(P^+_l, P^-_l)$ where
$\psi_\lambda$, $\psi_\mu$ vanish. Then we have that their ratio
$f:=\psi_\lambda/\psi_\mu$ still satisf\/ies the
conditions~\eqref{balf} at the remaining $m-r$ pairs of points and
has degree $\le m-r$ due to the cancelation of the zeros in the
denominator and numerator of $f$. Thus, the previous argument
applies and gives that $f$ is even. Therefore,
$\widetilde\psi/\psi$ is even.
\end{proof}

\subsection{Continuous limit}

As remarked in \cite{vD1}, the operator \eqref{l} with the
coupling parameters \eqref{co} in the continuous limit $\gamma\to
0$ turns into the $BC_1$-version \eqref{h} of the Lam\'e operator:
\begin{equation}\label{lim}
L=\mathrm{const} + \gamma^2 w(z)\circ H\circ w^{-1}(z) +
o(\gamma^2),\qquad \text{where}\quad w(z)=\prod_{p=0}^3
\left(\sigma_p(z)\right)^{g_p},
\end{equation}
and the coupling parameters $g_p$ are given by $g_p:=m_p+m'_p$.

To formulate a Bethe ansatz for the operator \eqref{h}, we put
$m=\sum\limits_{p=0}^3g_p$ and let $\psi(z)=\psi(z;k,t_1,\dots,
t_m)$ be the function \eqref{psi}. Let us impose the following $m$
relations on the parameters $k, t_1, \dots, t_m$:
\begin{equation}\label{cba}
\left[\frac{d^{2j-1}}{dz^{2j-1}}\left(\psi(z)e^{-m\eta_sz}\right)\right]_{z=\omega_s}=0
\qquad\text{for}\quad j=1,\dots, g_s\quad \text{and}\quad
s=0,\dots, 3.
\end{equation}

\begin{theorem}\label{cmt} Suppose the parameters $t_1,\dots, t_m, k$ satisfy the
Bethe equations \eqref{cba} and the conditions $t_i+t_j\notin
2\omega_1\Z+2\omega_2\Z$ for $1\le i\ne j\le m$. Then the function
$w^{-1}(z)\psi(z)$ given by \eqref{psi} with $w$ as in
\eqref{lim}, is an eigenfunction of the operator \eqref{h}.
\end{theorem}
This theorem is proved analogously to Theorem \ref{mt}.

\begin{example}
Let $g_p=0$ for $p=1,2,3$ and $g_0=m$. Then the operator \eqref{h}
becomes the Lam\'e operator $-d^2/dz^2+m(m+1)\wp(z)$. Its
eigenfunctions have the form $\psi(z)\sigma^{-m}(z)$ with
$\psi(z)=e^{kz}\prod\limits_{j=1}^m\sigma(z+t_j)$. The Bethe
ansatz equations \eqref{cba} for $k, t_1, \dots, t_m$ in this case
reduce to:
\begin{equation}\label{cbal}
\frac{d^{2j-1}\psi}{dz^{2j-1}}(0)=0 \qquad\text{for}\quad
j=1,\dots, m.
\end{equation}
We should note that this particular form of Bethe equations
dif\/fers from the classical result by Hermite~\cite{WW}. For
instance, in Hermite's equations one discards the points with
$t_i=t_j\mod 2\Gamma$, while in Theorem \ref{cmt} we discard the
points with $t_i=-t_j \mod 2\Gamma$. Thus, comparing Theorem
\ref{cmt} with the Hermite's result, we conclude that \eqref{cbal}
must be equivalent to Hermite's equations \cite{WW} provided
$t_i\pm t_j \notin 2\Gamma$ for $i\ne j$. The same remark applies
to the equations \eqref{cba} when compared to the Bethe ansatz in,
e.g., \cite{T}.
\end{example}

\subsection{Spectral curve}
Let us say few words about the structure of the solution set
$X\subset \c^m\times\c$ to the Bethe equations
\eqref{ba}--\eqref{ba'}. We will skip the details, since the
considerations here are parallel to those in \cite{FV1,FV2,KZ,Z}.

First, using the properties of $\sigma(z)$, one observes that
$\psi(z)$ acquires a constant factor under the transformations
\begin{equation*}\label{trans}
(t_1, \dots, t_m, k)\mapsto (t_1, \dots, t_j+2\omega_s, \dots,
t_m, k-2\eta_s)\qquad (s=1,2).
\end{equation*}
As a result, $X$ is invariant under these transformations. Also,
multiplying $\psi(z)$ by $e^{\pi i z/\gamma}$ does not af\/fect
the Bethe equations, because such an exponential factor is
(anti)periodic under the shifting of $z$ by multiples of $\gamma$.
Therefore, $X$ is invariant under the shifts of $k$ by $\pi
i/\gamma$:
\begin{equation*}\label{ktr}
(t_1, \dots, t_m, k)\mapsto (t_1, \dots, t_m, k+\pi i/\gamma).
\end{equation*}
Finally, $\psi$ does not change under permutation of $t_1,\dots,
t_m$, so $X$ is invariant under such permutations.

Let $\widetilde X$ denote the quotient of $X$ by the group
generated by all of the above transformations. Explicitly, let
$b_{s,j}(t_1, \dots, t_m)$ and $b'_{s,j}(t_1, \dots, t_m)$ denote
the left-hand side of equations \eqref{bae} and~\eqref{bae'}.
(Here $s=0,\dots, 3$ and $j=1,\dots, m_s$ or $j=1,\dots, m'_s$,
respectively.) Introduce the variable $q:=e^{2\gamma k}$. Then $X$
is described by the equations
\begin{equation}\label{xe}
b_{s,j}(t_1, \dots, t_m)=q^{2j}, \qquad b'_{s,j}(t_1, \dots,
t_m)=q^{2j-1}.
\end{equation}
Excluding the $q$-variable from the equations \eqref{xe}, we may
think of $\widetilde X$ as an algebraic subvariety in the
symmetric product $S^m\mathcal E$ of $m$ copies of the elliptic
curve $\mathcal E=\c/2\Gamma$ where
$\Gamma=\Z\omega_1+\Z\omega_2$. (See, however, Remark \ref{ree}
below.) Counting the number of equations, we conclude that every
irreducible component of $\widetilde X$ has dimension $\ge 1$.
Since we are interested (cf.\ Theorem~\ref{mt}) in those points
$(t_1,\dots, t_m)$ of $\widetilde X$ where $t_i+t_j\notin
2\Gamma$, we should restrict ourselves to the open part
$Y\subset\widetilde X$, lying in
\begin{equation*}
S^m\mathcal E\setminus \cup_{i<j}\{ t_i+t_j\equiv 0\
\mathrm{mod}\, 2\Gamma\}.
\end{equation*}
We need to show that $Y$ is nonempty and one-dimensional. To this
end, one easily observes from the equations
\eqref{bae}--\eqref{bae'} that the closure $\overline Y$ of $Y$ in
$S^m\mathcal E$ contains the points $P^+=(P_1^+, \dots, P_m^+)$
and $P^-=(P_1^-, \dots, P_m^-)$, in the notations of the proof of
the Lemma \ref{l2}. These `inf\/inite' points correspond to $q\to
0,\,\infty$ in \eqref{xe}. Similarly to \cite{FV2}, lemma 3.2, one
shows that near $P^\pm$ the variety $\overline Y$ looks like a
smooth curve, with $q^{\pm 1}$ being a local parameter. One can
show that~$\overline Y$ is an irreducible, projective curve, and
it should be regarded as the `spectral curve' for the operator
$L$. For every $(t_1,\dots, t_m)\in \overline Y\setminus\{P^+,
P^-\}$, the corresponding value of $q=e^{2\gamma k}$ is determined
from \eqref{xe}, and the corresponding $\psi(z)$ is unique, up to
a factor of the form $e^{\pi iNz/\gamma}$. We have
$L\psi=\epsilon\psi$, with the eigenvalue $\epsilon$ being a
single-valued function on $\overline Y$ which has two simple poles
at $P^\pm$. There is an involution $\nu$ on $\overline Y$, which
sends $(t_1,\dots, t_m)$ to $(-t_1, \dots, -t_m)$ and the
corresponding $\psi(z)$ to $\psi(-z)$; note that $\nu(P^+)=P^-$.
The function $\epsilon$ is $\nu$-invariant, and takes each its
value exactly twice on $\overline Y$. It is straightforward to
compute the asymptotics of $\epsilon$ and $\psi$ near $P^\pm$.
Finally, for generic value of $\epsilon$, the eigenspace of
meromorphic functions $\{f: Lf=\epsilon f\}$ is spanned by the
corresponding $\psi(z)$, $\psi(-z)$ over the f\/ield $K$ of
$2\gamma$-periodic meromorphic functions of $z$.

\begin{remark}\label{ree}
Note that in the case when all $m'_s=0$, the second set
\eqref{ba'} of the Bethe equations is absent, thus a shift
$k\mapsto k+\frac{\pi i}{2\gamma}$ is also allowed. In that case
the coef\/f\/icient $c(z)$ \eqref{c} vanishes, so $L$ has two
terms only, and it is easy to see that the transformation
$\psi\mapsto e^{\frac{\pi iz}{2\gamma}}\psi$ changes the sign of
the eigenvalue $\epsilon$. As a result, the subvariety of
$S^m\mathcal E$ which was obtained by excluding $q$ from
\eqref{xe}, will be a quotient of $\widetilde X$ by $\Z_2$, rather
than $\widetilde X$ itself (cf.~\cite{FV1,FV2,KZ,Z}).
\end{remark}

\pdfbookmark[1]{References}{ref}
\LastPageEnding

\end{document}